\documentclass{article}
\usepackage[utf8]{inputenc}
\usepackage{amsmath}
\usepackage{indentfirst}
\usepackage{amsthm}
\usepackage{tikz}
\usepackage{changepage}
\usepackage{url}
\usepackage{titling}
\usepackage{ragged2e}

\usepackage[maxbibnames=9,sortcites, doi=false, url=false, backend=biber]{biblatex}
 \usepackage{stackengine,amssymb,graphicx,scalerel}
\theoremstyle{definition}
\newtheorem{definition}{Definition}[section]
\numberwithin{definition}{section}
\newtheorem{theorem}{Theorem}[section]
\numberwithin{theorem}{section}

\newtheorem{lemma}{Lemma}[section]
\numberwithin{equation}{section}
\usetikzlibrary{backgrounds}
\theoremstyle{remark}
\newtheorem*{remark}{{Remark}}

\title{\bf Kinetic interaction Morawetz and concentration estimates}
\author{Nima Moini \\ University of California, Berkeley}
\date{}

\begin{document}
\begin{titlingpage}
\maketitle
\begin{abstract}
In this paper we will develop the kinetic interaction Morawetz estimate for mesoscopic evolutions. Furthermore, we will find a new a-priori bound and introduce the notion of interaction uncertainty. We will demonstrate that interaction uncertainty goes to infinity with time and that the interactions, averaged over time, concentrate within a specific collection of blind cones with an arbitrarily small apex angle. These results are solely based on conservation laws and thereby are all true for the Boltzmann equation and its different variants.
\end{abstract}
\setcounter{tocdepth}{1}
\tableofcontents
\end{titlingpage}

\section{Introduction}
In this paper we will develop new estimates associated to the evolution of interacting moving bodies on the mesoscopic scale. Consider a partial differential equation of the form:   
\begin{equation}
\begin{aligned}
    \partial_t f(x, \xi, t)+\xi.\nabla_x f(x, \xi, t) &= I(f, x, \xi, t) \\  f(x, \xi, 0)&=f_{0}(x,\xi)
\end{aligned}
\end{equation}
In the equation defined above, $x,\xi \in {\rm I\!R}^n$ and $t \in [0,\infty)$. We expect $f(x,\xi,t)$ to be non-negative and interpret it as the density or amplitude of particles with velocity $\xi$ located at $x$. The interaction $I(f,x,\xi,t)$ is expected to be mesoscopic in the sense defined below:

\begin{definition}\label{mesoscopic}
We say $I$ is a mesoscopic interaction if: 
\begin{align*}
                        \int_{{\rm I\!R}^n} I(f, x, \xi, t) \ d\xi&= 0 \\ \int_{{\rm I\!R}^n} I(f, x, \xi, t) |\xi|^2 \ d\xi&=0 \\ \int_{{\rm I\!R}^n} I(f, x, \xi, t) \xi \ d\xi&=0
\end{align*}
\end{definition}
\indent The characteristic of the partial differential equation (1.1) initiating from arbitrary $x,\xi \in {\rm I\!R}^n$ is the line $(x + t\xi,\xi, t) \subset {\rm I\!R}^n\times{\rm I\!R}^n\times {\rm I\!R}$. It is possible to interpret this equation as an infinite dimensional system of ordinary differential equations which start from the same initial value $f_{0}(x,\xi)$ and evolve along these characteristics: 
\begin{equation}\label{odeform}
\begin{aligned}
    \frac{d}{dt}f(x+t\xi,\xi,t)&=I(f,x+t\xi,\xi,t) \\
    f(x,\xi,0)&=f_{0}(x,\xi)
\end{aligned}
\end{equation}
We will assume that at time zero the total mass, energy and momentum defined below are finite:
\begin{equation}\label{mme}
\begin{aligned}
  ({\it mass})&  &{\bf M} &=\int_{{\rm I\!R}^n}\int_{{\rm I\!R}^n} f(x,\xi, 0) \ dxd\xi < \ \infty  \\
     ({\it energy})&  & {\bf E} &=\int_{{\rm I\!R}^n}\int_{{\rm I\!R}^n} f(x,\xi, 0)|\xi|^2 \ dxd\xi < \ \infty
     \\
     ({\it momentum})&  &{\bf V} &=\int_{{\rm I\!R}^n}\int_{{\rm I\!R}^n} f(x,\xi, 0)\xi \ dxd\xi < \ \infty
\end{aligned}
\end{equation}
\noindent For a subset of the following results we will also assume a localization of mass at time zero in the sense defined below: 
\begin{align}\label{localization}
    \int_{{\rm I\!R}^n}\int_{{\rm I\!R}^n}f(x, \xi, 0)|x|^2 dxd\xi < \infty
\end{align}

\indent One immediate consequence of the equation above and the definition of a mesoscopic interaction is that, for the solutions of equation (1.1), the total amount of mass, momentum and energy remain invariant. For example, the conservation of mass can be proven by the argument below:
\begin{multline*}
    \frac{d}{dt}\int_{{\rm I\!R}^n}\int_{{\rm I\!R}^n} f(x, \xi, t)  \ dxd\xi= 
\int_{{\rm I\!R}^n}\int_{{\rm I\!R}^n}\frac{d}{dt}f(x+t\xi,\xi,t)\ dxd\xi=\\ \int_{{\rm I\!R}^n}\int_{{\rm I\!R}^n} I(f, x+t\xi, \xi, t)  \ dxd\xi=\int_{{\rm I\!R}^n}\int_{{\rm I\!R}^n} I(f, x, \xi, t)  \ dxd\xi=0 
\end{multline*}

The Boltzmann equation has a similar structure to the equation (1.1) for a specific interaction $I$. The notation $Q(f,f)(x,\xi,t)$ represents the interactions of the particles and is called the Boltzmann collision operator \cite{MR1313028, MR1307620}. In the case which this operator represents collisions of hard elastic spheres it obtains the form below:
\begin{gather*}
    Q(f, f)(x,\xi,t) = \int_{{\rm I\!R}^n}\int_{S^{n-1}}(f^{'} f^{'}_{*} - f f_{*})|{\bf n}.(\xi- \xi_{*})| d{\bf n}d\xi_*  
 \end{gather*}
In the expression above, ${\bf n}$ is the normal unit vector to  the $n-1$ dimensional unit sphere $S^{n-1}$ and, as is customary in the field, we used the notation below to describe the velocity of particles before and after the collisions: 
\begin{gather*}
    \xi^{'} = \xi - {\bf n}.(\xi- \xi_{*}){\bf n}  \\
     \xi^{'}_{*} = \xi_{*} + {\bf n}.(\xi- \xi_{*}){\bf n}  \\ 
f = f(x,\xi,t),\ f_{*} = f(x,\xi_*,t) \\ f^{'} = f(x, \xi^{'}, t) , \ f^{'}_{*} = f(x, \xi^{'}_{*}, t) 
\end{gather*}
\noindent Assume that $\xi$ and $\xi_{*}$ are velocities of two colliding elastic balls with unit mass, and let ${\bf n}$ or equivalently $-{\bf n}$, represent the unit normal vector to the plane which uniquely describes the relative position of these two spheres upon the collision. Then $\xi^{'}$ and $\xi_{*}^{'}$ are the velocities of the  particles after the collision. For a fixed ${\bf n}$, these velocities are unique solutions to the conservation laws of momentum and energy written below: 
\begin{align*}
    |\xi|^2 + |\xi_{*}|^2 &= |\xi^{'}|^2 + |\xi^{'}_{*}|^2 \\ \xi + \xi_{*} &= \xi^{'}+\xi_{*}^{'}
\end{align*}

As shown by Boltzmann, the collision operator satisfies the conditions of Definition \ref{mesoscopic} and therefore is a mesoscopic interaction. We will replicate his argument here. Consider the  well-known change of variables below:
\begin{align*}
    (\xi, \xi_{*}) \Longrightarrow (\xi_* , \xi), \ \ \ (\xi, \xi_{*}) \Longrightarrow (\xi^{'} , \xi_{*}^{'}), \ \ \   (\xi, \xi_{*}) \Longrightarrow (\xi^{'}_* , \xi^{'})  
\end{align*}
The transformations above are measure preserving. They represent the intrinsic symmetries of the conservation laws, in which the role between particles as well as the velocity of particles before and after a collision are indistinguishable. We will implement these change of variables for the Boltzmann collision operator and an arbitrary $\phi(\xi)$, we get: 
\begin{multline*}
  \int_{{\rm I\!R}^n} Q(f,f)(x,\xi,t)\phi(\xi) \ d\xi =  \frac{1}{4}\int_{{\rm I\!R}^n}\int_{{\rm I\!R}^n}\int_{S^{n-1}}(f^{'} f^{'}_{*} - f f_{*}) \\ \times \big(\phi(\xi) + \phi(\xi_{*}) - \phi(\xi^{'})  - \phi(\xi^{'}_{*})\big)|{\bf n}.(\xi- \xi_{*})| \ d{\bf n} d\xi_* d\xi
\end{multline*}
Set $\phi(\xi)$ equal to either 1, $|\xi|^2$ or $\xi_{i}$ for $1\leq i \leq n$. It follows from the conservation laws that $\phi(\xi) + \phi(\xi_{*}) - \phi(\xi^{'})  - \phi(\xi^{'}_{*})=0$. We conclude that assumptions of Definition \ref{mesoscopic} are satisfied by $Q$ and hence it is a mesoscopic interaction. \\

The existence of global classical solutions to the Boltzmann equation for small data has been shown before \cite{MR760333}. Additionally, there exist a scattering theory for small data with convergence in the $L^{\infty}$ setting and asymptotic completeness results \cite{moini2021evolution}. A concept of renormalized solutions exists that provides a general existence theory of weak solutions for the Boltzmann equation \cite{MR1014927}. The counterpart of this renormalization in the microscopic scale is unknown. This missing link makes it infeasible to completely relate these weak solutions to the underlying physical phenomena which the Boltzmann equation is intended to describe. The existence of global in time classical solutions to the Boltzmann equation for the large data or examples of finite time blow up are still important open problems.\\

\indent The results of this paper are solely based on the conservation laws and are independent of the specific structure of the interaction, therefore they are also true for the Boltzmann equation. The kinetic descriptions of conservation laws has been studied before \cite{ MR2083859,MR1201239,MR1284790,MR1266203,MR2397052} using methods like velocity averaging lemmas and the notion of entropy \cite{MR1127927,MR923047,MR2043752}. In this paper we will not use these conventional methods. We will generalize some of the recent results in the kinetic theory from single particle estimates to the interaction or 2-particle estimates. We will start with a short summery of the existing results appearing in \cite{moini2021evolution}. Consider the definitions below: 
\begin{equation}\label{singleparticle}
\begin{aligned}
    (angular \ momentum)&& A(t) &= \int_{{\rm I\!R}^n}\int_{{\rm I\!R}^n} f(x,\xi, t)x.\xi \ dxd\xi \\ 
    (uncertainty)&& U(t) &= \int_{{\rm I\!R}^n}\int_{{\rm I\!R}^n} f(x,\xi, t)|x||\xi| \ dxd\xi \\ 
    (relative \ angular \ norm)&& \|f\|_{G}&=\sup_{t} \big(U(t) - A(t)\big)
\end{aligned}
\end{equation}
For any positive solution of the equation (1.1) subject to a mesoscopic interaction, angular momentum increases linearly:   
\begin{align*}
    A(t) = A(0) + t{\bf E}
\end{align*}
Furthermore,  we have: 
\begin{align}\label{ag}
    \|f\|_G \leq \int_{{\rm I\!R}^n}\int_{{\rm I\!R}^n} f(x, \xi, 0)|x|^2 dxd\xi + {\bf E} - A(0)
\end{align}
Which implies, the uncertainty increases on average linearly proportional to the amount of the total energy and goes to infinity with time:
\begin{align*}
        \lim_{t\rightarrow\infty}U(t) = \infty 
\end{align*}
\indent The uncertainty defined above has a physical interpretation with similarities to its quantum counterpart  \cite{moini2021evolution}. Theses results led to a new physical intuition based on the notion of a blind cone with respect an observer. These blind cones have been used in \cite{moini2021evolution} to prove that the total energy within any bounded set of the spatial variable is integrable over time. Assume $D \subset {\rm I\!R}^n$ is a bounded set of the spatial variable, we have:
\begin{align}\label{mes_morawetz}
 ({\it Kinetic \ Morawetz \ estimate}) \ \int_{0}^{\infty}\int_{{\rm I\!R}^n}\int_{D} f(x,\xi, t)|\xi|^2 \ dxd\xi dt < \infty
\end{align}

\indent The previous result has analogies to the Morawetz estimate for the non linear Schrodinger's equation \cite{MR2233925}. The convergence of the time integral in (\ref{mes_morawetz}) implies that the total mass of  particles with velocity greater than any positive fixed positive number like $v$ is integrable over time as well. Including more assumptions will further illuminate the physical intuition behind the estimate. Assume $f$ is bounded and has a bounded derivative with respect to the time. The previous argument indicates that inside any bounded set of the spatial variable like $D \subset {\rm I\!R}^n$, as time goes to infinity, almost every particle with a magnitude of velocity greater than $v$ will inevitably leave the bounded set. Thus, as time goes to infinity, an idle observer will almost only identify particles with velocity $\xi = 0$ inside any bounded set. The Galilean invariance of the setting will generalize this argument to bounded moving regions and moving observers. Each particle has a tendency to travel with other particles whose velocities are indistinguishable from each other and as a result avoid interaction with particles that move at different velocities. \\

\indent In {\bf Section 2}, we will generalize (\ref{mes_morawetz}) and create the {\it 2-particle} or {\it interaction  Morawetz estimate} (Theorem \ref{kime}) for the evolution of mesoscopic interactions. This estimate has analogies to the interaction Morawetz estimate for non linear Schrodinger's equation \cite{doi:10.1137/19M1270586, 10.2307/40345365,MR2233925,MR3229599} and has physical interpretations. This estimate implies that, along the evolution of solutions to the equation (1.1) and as the time increases, only particles with identical velocities can remain within an arbitrary fixed distance of one another. \\ 

\indent In {\bf Section 3}, we will generalize the notions appearing in (1.5). We will demonstrate a new a-priori bound (Definition \ref{ian}) for the 2-particle interactions and introduce a notion of interaction uncertainty (Definition \ref{intunc}). We will show that the interaction uncertainty, similar to the single particle case, will go to infinity with time (Theorem \ref{ianthm}). These results will be used in Section 4.\\

 \indent It is possible to use the kinetic Morawetz estimate (\ref{mes_morawetz}) together with the bound (\ref{ag}) for the relative angular norm to prove concentration of mass overtime in a specific sense \cite{moini2021evolution}. This concentration implies that the particles will eventually move away from any observer in a radial manner. Therefore, the angular part of the gradient will vanish relative to any observer. Consequently, the total mass concentrates, in an averaged sense over time, within the collection of blind cones with respect to any fixed observer.  \\
 
 \indent In {\bf Section 4}, we will create an analogy (Theorem \ref{conc}) to the described concentration of mass for the case of interactions. We will show that the interactions, averaged over time, will concentrate within a specific collection of the blind cones (Definition \ref{gamma}) with arbitrary small apex angle. The relation between Morawetz and interaction Morawetz estimates bears resemblance to the relationship between concentration of mass and concentration of interactions. The concentration of interactions within the collection of blind cones implies that, averaged over time, almost every collision occurs between particles with arbitrary close velocities or arbitrary small angle of deflection. These type of interactions are often called the grazing collisions.

\subsection*{Acknowledgements} I am grateful for the mentorship of Daniel Tataru during the creation of this work. I want to express my gratitude for useful discussions with Fraydoun Rezakhanlou and Mehrdad Shahshahani. This paper would not have been written without the love and support of Amanda Jeanne Tose (AJT), Kaveh and Zahra.

\newpage

\section{Interaction Morawetz estimate}
In this section we will develop a kinetic interaction Morawetz estimate based on the conservation laws of the classical physics. 

\begin{definition} Let $A_{L}(t)$ be the {\it 2-particle localized angular momentum} associated to  $f(x, \xi, t)$ at time $t$:
\begin{align*}
    A_{L}(t) = \iiiint\limits_{{\rm I\!R}^{4n}} f(x, \xi, t)f(x_{0}, \xi_{0}, t) \frac{x - x_{0}}{|x - x_0|} \cdot (\xi - \xi_{0}) \ dx d\xi dx_{0} d\xi_{0}
\end{align*}
\end{definition}
\begin{remark}
The quantity defined above is bounded over time. 
\begin{align}\label{albound}
A_{L}(t) < ({\bf M} + {\bf E})^2
\end{align}
\noindent It is possible to interpret $A_{L}(t)$ as a measure of dispersion for the pairs of particles and has been used in different contexts \cite{MR2229996, MR1345903}. The single particle variant of this quantity appears in \cite{moini2021evolution} and justifies the language used above. 
\end{remark}

\begin{theorem}\label{kime}[{\it Kinetic interaction Morawetz estimate}] Assume $f$ is a non negative solution of the equation (1.1) subject to a mesoscopic interaction. Furthermore, assume the total mass ({\bf M}) and energy ({\bf E}) are as in (\ref{mme}). There exist a constant $W$ depending only on the total mass and energy such that for any positive $R$ we have: 
\begin{multline*}
   \frac{1}{R}\int_{0}^{\infty} \int_{{\rm I\!R}^{n}}\int_{{\rm I\!R}^{n}}\int_{B(x_0, R)}\int_{{\rm I\!R}^{n}} f(x, \xi, t)f(x_{0}, \xi_{0}, t) |\xi - \xi_0|^2 \ d\xi dx d\xi_{0} dx_0 dt  \\ < W({\bf M}, {\bf E})
   \end{multline*}
\begin{proof}
We start by differentiating $A_L(t)$ with respect to time:
\begin{equation}\label{dadt}
\begin{aligned}
 \frac{d}{dt} A_L = \iiiint\limits_{{\rm I\!R}^{4n}} \partial_{t}f(x, \xi, t) \frac{x - x_{0}}{|x - x_0|} \cdot (\xi - \xi_{0}) \ dxd\xi \  f(x_{0}, \xi_{0}, t) \ dx_{0} d\xi_{0}\\+\iiiint\limits_{{\rm I\!R}^{4n}} \partial_{t}f(x_{0}, \xi_{0}, t) \frac{x - x_{0}}{|x - x_0|} \cdot (\xi - \xi_{0}) \ dx_{0}d\xi_{0} \  f(x, \xi, t) \ dx d\xi   
\end{aligned}
\end{equation}
For the inner double integral appearing in the first term  of the right hand side of the equation above we have: 
\begin{multline*}
     \int_{{\rm I\!R}^n}\int_{{\rm I\!R}^n} \partial_t f(x,\xi,t)\frac{x - x_{0}}{|x - x_0|} \cdot (\xi - \xi_{0}) \ dxd\xi\\ + \int_{{\rm I\!R}^n}\int_{{\rm I\!R}^n}\frac{x - x_{0}}{|x - x_0|} \cdot (\xi - \xi_{0}) \xi.\nabla_{x}f(x,\xi, t)\ dxd\xi \\= \int_{{\rm I\!R}^n}\int_{{\rm I\!R}^n} I(f, x,\xi,t)\frac{x - x_{0}}{|x - x_0|} \cdot (\xi - \xi_{0}) \ dxd\xi=0
\end{multline*}
Therefore:   
\begin{multline*}
     \int_{{\rm I\!R}^n}\int_{{\rm I\!R}^n} \partial_t f(x, \xi, t)\frac{x - x_{0}}{|x - x_0|} \cdot (\xi - \xi_{0}) dxd\xi \\ = - \int_{{\rm I\!R}^n}\int_{{\rm I\!R}^n}\frac{x - x_{0}}{|x - x_0|} \cdot (\xi - \xi_{0}) \xi.\nabla_{x}f(x,\xi, t)\ dxd\xi
\end{multline*}
Similarly, for the inner double integral appearing in the second term of the equation (\ref{dadt}) we get: 
\begin{multline*}
         \int_{{\rm I\!R}^n}\int_{{\rm I\!R}^n} \partial_t f(x_0, \xi_0, t) \frac{x - x_{0}}{|x - x_0|} \cdot (\xi - \xi_{0}) dx_0 d\xi_0 \\ = - \int_{{\rm I\!R}^n}\int_{{\rm I\!R}^n} \frac{x - x_{0}}{|x - x_0|} \cdot (\xi - \xi_{0}) \xi_0.\nabla_{x}f(x_0, \xi_0, t) \ dx_0d\xi_0
\end{multline*}
We will continue with integration by parts with respect to $x$ for the two former equations:  
\begin{multline*} 
    \partial_t \int_{{\rm I\!R}^n}\int_{{\rm I\!R}^n}   f(x,\xi,t)\frac{x - x_{0}}{|x - x_0|} \cdot (\xi - \xi_{0})\ dxd\xi\\= \int_{{\rm I\!R}^n}\int_{{\rm I\!R}^n} \frac{\xi \cdot (\xi - \xi_0)|x - x_0|^2 - ((x- x_0) \cdot \xi) ((x - x_0)\cdot (\xi - \xi_0))}{|x-x_0|^3} \\ \times f(x,\xi, t) \ dx d\xi
\end{multline*}
Likewise: 
\begin{multline*} 
    \partial_t \int_{{\rm I\!R}^n}\int_{{\rm I\!R}^n}   f(x_0,\xi_0,t)\frac{x - x_{0}}{|x - x_0|} \cdot (\xi - \xi_{0})\ dx_0 d\xi_0 \\ = \int_{{\rm I\!R}^n}\int_{{\rm I\!R}^n} \frac{-\xi_0 \cdot (\xi - \xi_0)|x - x_0|^2 + ((x- x_0) \cdot \xi_0) ((x - x_0)\cdot (\xi - \xi_0))}{|x-x_0|^3} \\ \times f(x_0,\xi_0, t) \ dx_0 d\xi_0
\end{multline*}
Plug in the the computations above in the equation (\ref{dadt}). We will obtain the expression below, where $\theta(x, \xi)$ is the angle between the two vectors:
\begin{multline*}
        \frac{d}{dt} A(t) = \iiiint\limits_{{\rm I\!R}^{4n}} f(x, \xi, t)f(x_{0}, \xi_{0}, t) \frac{1}{|x - x_0|}|\xi - \xi_0|^2 \\ \times \sin^2{(\theta(x-x_0, \xi-\xi_0)}) \ d\xi dx d\xi_{0} dx_0
\end{multline*}
This shows that $A_{L}(t)$ has a positive derivative and is monotone. Since this quantity is bounded as (\ref{albound}), we get: 
\begin{multline*}
    \int_{0}^{\infty} \iiiint\limits_{{\rm I\!R}^{4n}} f(x, \xi, t)f(x_{0}, \xi_{0}, t) \frac{1}{|x - x_0|}|\xi - \xi_0|^2 \\ \times \sin^2{(\theta(x-x_0, \xi-\xi_0)}) \ d\xi dx d\xi_{0} dx_0 dt < ({\bf M}+{\bf E})^2
\end{multline*}
We will continue by defining a particular kind of cone. Let $C_{x_0, \xi_0}(x, c) \subset {{\rm I\!R}^{n}}$ be the blind cone with apex angle $c>0$ at point $x \in {{\rm I\!R}^{n}}$ with respect to a moving observer at $x_0 \in {{\rm I\!R}^{n}}$  with velocity $\xi_0 \in {{\rm I\!R}^{n}}$:
\begin{equation}\label{bcone}
\begin{aligned}
       C_{x_0, \xi_0}(x, c) = \{\xi \in {\rm I\!R}^n \big | \theta(x-x_0, \xi - \xi_0) \notin  [c, \pi-c]\}  
\end{aligned}
\end{equation}
We will identify this cone as a subset of the space of velocities at point $x$. By removing the blind cones $C_{x_0, \xi_0}(x, c)$ from the space of velocities for every point $x$ and using the previous computations we get that for any positive $R$:
\begin{equation}\label{blremoved}
\begin{aligned}
    \frac{\sin^2(c)}{R}\int_{0}^{\infty} \int_{{\rm I\!R}^{n}}\int_{{\rm I\!R}^{n}}\int_{B(x_0, R)}\int_{{\rm I\!R}^{n} - C_{x_0, \xi_{0}(x, c)}} f(x, \xi, t)f(x_{0}, \xi_{0}, t) \\ \times |\xi - \xi_0|^2 \ d\xi dx d\xi_{0} dx_0 dt < ({\bf M}+{\bf E})^2
\end{aligned}
\end{equation}

\noindent Now choose any three distinct {\it observers}: $O_{1}, O_{2}, O_{3} \in \partial B(0, R)$ and for an arbitrary $x\in B(0, R)$ define $P=\partial B(0, R)\cap C_{O_1, \xi_{0}}(x, 2c)$. If there exists an $O_i$ such that $O_{i} \notin P$ then the blind cones $C_{O_i, \xi_0}(x, c)$ and $C_{O_1, \xi_0}(x, c)$ have an empty intersection. Set $P$ belongs to $\partial B(0, R)$ and is made of two path connected components, consider the longest short path on each component and set $K$ to be the maximum length of the two. For any fixed $R$ it is possible to choose $c$ small enough such that $K$ becomes as small as desired. Now set $c$ small enough such that $K$ becomes smaller than the shortest path on the sphere between any two of the there observers. The pigeon hole principle implies that, since each path connected component of $P$ can only contain maximum one of the observers, there exists an $O_{i}$ such that $O_{i} \notin P$. The previous argument implies that for any $x \in B(0, R)$ the blind cones with respect to the three observers have an empty intersection:   
\begin{align*}
    C_{O_1, \xi_0}( x, c) \cap  C_{O_2, \xi_0}( x, c) \cap  C_{O_3, \xi_0}(x, c) = \emptyset      
\end{align*}
Therefore, for any $ x_0 \in {\rm I\!R}^{n}$ we have:
\begin{equation}
\begin{aligned}\label{3observers}
    C_{x_0 + O_1, \xi_0}( x, c) \cap  C_{x_0 + O_2, \xi_0}( x, c) \cap  C_{x_0 + O_3, \xi_0}( x, c) = \emptyset      
\end{aligned}
\end{equation}
Consider the definitions below for $J_1$, $J_2$ and $J_3$. The following bounds are consequences of (\ref{blremoved}):
\begin{multline*}
    J_1 = \frac{\sin^2(c)}{R}\int_{0}^{\infty}\int_{{\rm I\!R}^{n}}\int_{{\rm I\!R}^{n}}\int_{B(x_0+O_1, 2R)}\int_{{\rm I\!R}^{n} - C_{x_0+O_1, \xi_{0}}(x, c)} f(x, \xi, t)\\ \times f(x_{0}, \xi_{0}, t)|\xi - \xi_0|^2 \ d\xi dx d\xi_{0} dx_0 dt < ({\bf M}+ {\bf E})^2
\end{multline*}
\begin{multline*}    J_2 =  \frac{\sin^2(c)}{R}\int_{0}^{\infty}\int_{{\rm I\!R}^{n}}\int_{{\rm I\!R}^{n}}\int_{B(x_0+O_2, 2R)}\int_{{\rm I\!R}^{n} - C_{x_0+O_2, \xi_{0}}(x, c)} f(x, \xi, t)\\ \times f(x_{0}, \xi_{0}, t) |\xi - \xi_0|^2 \ d\xi dx d\xi_{0} dx_0 dt < ({\bf M}+ {\bf E})^2
\end{multline*}
\begin{multline*}
    J_3 =  \frac{\sin^2(c)}{R}\int_{0}^{\infty}\int_{{\rm I\!R}^{n}}\int_{{\rm I\!R}^{n}}\int_{B(x_0+O_3, 2R)}\int_{{\rm I\!R}^{n} - C_{x_0+O_3, \xi_{0}}(x, c)} f(x, \xi, t)\\ \times f(x_{0}, \xi_{0}, t) |\xi - \xi_0|^2 \ d\xi dx d\xi_{0} dx_0 dt < ({\bf M}+ {\bf E})^2
\end{multline*}

\noindent From (\ref{3observers}) we know that it is possible to choose the apex angle $c$ small enough such that for any point $x \in B(x_0, R)$ the blind cones with respect to the three observers $x_0+O_1$, $x_0+O_2$ and $x_0+O_3$ intersect trivially. Also consider that:   
\begin{align*}
    B(x_0, R) \subset B(x_0 + O_1, 2R) \cap B(x_0 + O_2, 2R) \cap B(x_0 + O_3, 2R)
\end{align*}

\noindent This shows that any subset of $B(x_0, R) \times {\rm I\!R}^n$ for all $x_0 \in {\rm I\!R}^n$ is covered at least  once in the domains of integration for $J_1, J_2$ and $J_3$. Therefore the positivity of integrands completes the proof, there exists some constant $W({\bf M}, {\bf E})$ such that: 
\begin{multline*}
   \frac{1}{R}\int_{0}^{\infty} \int_{{\rm I\!R}^{n}}\int_{{\rm I\!R}^{n}}\int_{B(x_0, R)}\int_{{\rm I\!R}^{n}} f(x, \xi, t)f(x_{0}, \xi_{0}, t) |\xi - \xi_0|^2 \ d\xi dx d\xi_{0} dx_0 dt \\ < J_1 + J_2 + J_3 < W({\bf M}, {\bf E})
\end{multline*}
\end{proof}
\end{theorem}
\begin{remark} An illustration of the blind cones with respect to an idle observer ($\xi_0$=0) at the origin ($x_0=0$) can be found in \cite{moini2021evolution}. 
\end{remark}
\newpage

\section{New a priori bound and the interaction uncertainty}
 In this section we will find a new a-priori bound for the evolution of the mesoscopic interactions. We will introduce the notion of uncertainty associated to the interactions and demonstrate that as time goes to infinity the interaction uncertainty goes to infinity as well. These results and the kinetic interaction Morawetz estimate (Theorem \ref{kime}) will be used in Section 4 to obtain a concentration result for the interactions. 
\begin{lemma}\label{lemmaunc}
Assume $f$ is a solution of the equation (1.1) subject to a mesoscopic interaction and that the assumptions (\ref{mme}) and (\ref{localization}) are true. We have: 
\begin{multline*}
   \iiiint\limits_{{\rm I\!R}^{4n}}f(x, \xi, t)f(x_{0}, \xi_{0}, t) \ (x - x_0) \cdot (\xi - \xi_{0}) \ dx d\xi dx_0 d\xi_0 \\ = \iiiint\limits_{{\rm I\!R}^{4n}}f(x, \xi, 0)f(x_{0}, \xi_{0}, 0) \ (x - x_0) \cdot (\xi - \xi_{0}) \ dx d\xi dx_0 d\xi_0 \\+ t\iiiint\limits_{{\rm I\!R}^{4n}}f(x, \xi, 0)f(x_{0}, \xi_{0}, 0) \ |\xi - \xi_0|^2 \ dx d\xi dx_0 d\xi_0
   \end{multline*}
\begin{proof} Start with the change of variables $x \rightarrow x+t\xi$ and $x_0 \rightarrow t\xi_0$: 
\begin{multline*} 
   \iiiint\limits_{{\rm I\!R}^{4n}}f(x, \xi, t)f(x_{0}, \xi_{0}, t) \ (x - x_0) \cdot (\xi - \xi_{0}) \ dx d\xi dx_0 d\xi_0 \\=\Big(\iiiint\limits_{{\rm I\!R}^{4n}}f(x+t\xi, \xi, t)f(x_{0}+t\xi_{0}, \xi_{0}, t) \ (x - x_0 ) \cdot (\xi - \xi_{0}) \ dx d\xi dx_0 d\xi_0 \\+ t \iiiint\limits_{{\rm I\!R}^{4n}}f(x+t\xi, \xi, t)f(x_{0}+t\xi_{0}, \xi_{0}, t) \ |\xi - \xi_0|^2 \ dx d\xi dx_0 d\xi_0\Big)
\end{multline*}
Continue with the computation below: 
\begin{multline*}
     \frac{d}{dt}\iiiint\limits_{{\rm I\!R}^{4n}}f(x, \xi, t)f(x_{0}, \xi_{0}, t) \ (x - x_0) \cdot (\xi - \xi_{0}) \ dx d\xi dx_0 d\xi_0 \\= 2\iiiint\limits_{{\rm I\!R}^{4n}} I(x+t\xi, \xi, t)f(x_{0}+t\xi_{0}, \xi_{0}, t) \ (x - x_0 ) \cdot (\xi - \xi_{0}) \ dx d\xi dx_0 d\xi_0 \\+ 2t \iiiint\limits_{{\rm I\!R}^{4n}}I(x+t\xi, \xi, t)f(x_{0}+t\xi_{0}, \xi_{0}, t) \ |\xi - \xi_0|^2 \ dx d\xi dx_0 d\xi_0 \\+ \iiiint\limits_{{\rm I\!R}^{4n}}f(x+t\xi, \xi, t)f(x_{0}+t\xi_{0}, \xi_{0}, t) \ |\xi - \xi_0|^2 \ dx d\xi dx_0 d\xi_0
\end{multline*}
We will compute each of the three terms appearing in the former equation. Because $I$ is a mesoscopic interaction (Definition \ref{mesoscopic}), for the first term we have:
\begin{multline*}
    2\iiiint\limits_{{\rm I\!R}^{4n}} I(x+t\xi, \xi, t)f(x_{0}+t\xi_{0}, \xi_{0}, t) \ (x - x_0 ) \cdot (\xi - \xi_{0}) \ dx d\xi dx_0 d\xi_0 \\ =  2\iiiint\limits_{{\rm I\!R}^{4n}} I(x, \xi, t)f(x_{0}, \xi_{0}, t) \ (x - x_0 + t(\xi_0 - \xi) \cdot (\xi - \xi_{0}))) d\xi dx dx_0 d\xi_0 =0
\end{multline*}
Similarly for the second term we get:
\begin{multline*}
    2\iiiint\limits_{{\rm I\!R}^{4n}} I(x+t\xi, \xi, t)f(x_{0}+t\xi_{0}, \xi_{0}, t) |\xi - \xi_0|^2 \ dx d\xi dx_0 d\xi_0 \\ =  2\iiiint\limits_{{\rm I\!R}^{4n}} I(x, \xi, t)f(x_{0}, \xi_{0}, t) \ |\xi - \xi_0|^2 d\xi dx dx_0 d\xi_0 =0
\end{multline*}
\noindent Finally we will compute the third term: 
\begin{multline*}
    \iiiint\limits_{{\rm I\!R}^{4n}}f(x+t\xi, \xi, t)f(x_{0}+t\xi_{0}, \xi_{0}, t) \ |\xi - \xi_0|^2 \ dx d\xi dx_0 d\xi_0 \\ =\iiiint\limits_{{\rm I\!R}^{4n}}f(x, \xi, 0)f(x_{0}, \xi_{0}, 0) \ |\xi - \xi_0|^2 \ dx d\xi dx_0 d\xi_0\\+ 2\iiiint\limits_{{\rm I\!R}^{4n}}\int_{0}^{t}I(x+s\xi, \xi, s)f(x_{0}, \xi_{0}, 0) \ |\xi - \xi_0|^2 \ ds dx d\xi dx_0 d\xi_0\\+\iiiint\limits_{{\rm I\!R}^{4n}}\int_{0}^{t}\int_{0}^{t}I(x+s\xi, \xi, s)I(x_{0}+z\xi_{0}, \xi_{0}, z) \ |\xi - \xi_0|^2 \ ds dz dx d\xi dx_0 d\xi_0
\end{multline*}
As a consequence of Definition \ref{mesoscopic}, after a change of variables and changing the order of integrations, we conclude that the last two terms of the right hand side of the equation above are zero. Therefore we complete the proof: 
\begin{multline*}
      \frac{d}{dt}\iiiint\limits_{{\rm I\!R}^{4n}}f(x, \xi, t)f(x_{0}, \xi_{0}, t) \ (x - x_0) \cdot (\xi - \xi_{0}) \ dx d\xi dx_0 d\xi_0 \\= \iiiint\limits_{{\rm I\!R}^{4n}}f(x, \xi, 0)f(x_{0}, \xi_{0}, 0) \ |\xi - \xi_0|^2 \ dx d\xi dx_0 d\xi_0
\end{multline*}
\end{proof}
\end{lemma}
\begin{definition}\label{intunc}
Assume $f$ is a non negative solution of the equation (1.1). Let $U_{I}(t)$ be the {\it interaction uncertainty} associated to $f$ at time $t$ defined as: 
\begin{align*}
       U_{I}(t)  = \iiiint\limits_{{\rm I\!R}^{4n}}f(x, \xi, t)f(x_{0}, \xi_{0},t)|x-x_0||\xi-\xi_0|dxd\xi dx_0d\xi_0
\end{align*}
\end{definition}
\begin{remark}
It is possible to interpret the interaction uncertainty defined above analogous to its quantum mechanics counterpart, that is by bringing attention to a fundamental limit on the precision of the physical measurements. Consider an idle observer observer located at the origin and assume the speed of light is $C$. For each particle located at $x$ there is a minimum delay of $T = C^{-1}|x|$ between the actual time of measurement and observation at the origin. The quantity $T \times f(0, 0, t)f(x, \xi, t)|x|=C^{-1}f(0, 0, t)f(x, \xi, t)|x||\xi|= $ represents the uncertainty of measurement relative to an idle observer at the origin, due to this
interval of delay. The role of the observer and the particle are symmetric in the context of interactions. Therefore after including all the particles and observers, we will get the definition above for $C=1$. 

\end{remark}
\begin{definition}\label{ian} Assume $f$ is a non negative solution of the equation (1.1). Let $\|\|_{IG}$ be the {\it interaction relative  angular norm} associated to $f$ defined as:
\begin{multline*}
   \|f\|_{IG}  = \sup_{t} \iiiint\limits_{{\rm I\!R}^{4n}}f(x, \xi, t)f(x_{0}, \xi_{0}, t) \Big(|x - x_0||\xi - \xi_0| \\- (x - x_0) \cdot (\xi - \xi_0)\Big) \ d\xi dx d\xi_{0} dx_0
\end{multline*}
\end{definition}
\begin{theorem}\label{ianthm} Assume $f$ is a non negative solution of equation (1.1). Additionally, assume that the assumptions (\ref{mme}) and (\ref{localization}) are true. We have: 
\begin{multline*}
   \|f\|_{IG} < \iiiint\limits_{{\rm I\!R}^{4n}}f(x, \xi, 0)f(x_{0}, \xi_{0}, 0) \Big(|x - x_0|^2+|\xi - \xi_0|^2 \\ - (x-x_0).(\xi - \xi_0)\Big)d\xi dx d\xi_0 dx_0
\end{multline*}
Furthermore, the interaction uncertainty goes to infinity with time:
\begin{align*}
    \lim_{t\rightarrow\infty}U_{I}(t) = \infty
\end{align*}
\begin{proof} Implement the change variables $x \rightarrow x + t\xi$ and $x_0 + t\xi_0$ followed by the triangle inequality, we get: 
\begin{multline*}
    \iiiint\limits_{{\rm I\!R}^{4n}}f(x, \xi, t)f(x_{0}, \xi_{0}, t) \Big(|x - x_0||\xi - \xi_0| - (x - x_0) \cdot (\xi - \xi_0)\Big) \ d\xi dx d\xi_{0} dx_0 \\ \leq \iiiint\limits_{{\rm I\!R}^{4n}}f(x+t\xi, \xi, t)f(x_{0}+t\xi_0, \xi_{0}, t) \Big(|x - x_0||\xi - \xi_0| + t|\xi - \xi_0)|^2 \\- (x - x_0+ t(\xi - \xi_0)) \cdot (\xi - \xi_0)\Big) \ d\xi dx d\xi_{0} dx_0
\end{multline*}
Therefore, using the non negativity of $f$ we have: 
\begin{multline*}
     \iiiint\limits_{{\rm I\!R}^{4n}}f(x, \xi, t)f(x_{0}, \xi_{0}, t)\Big(|x - x_0||\xi - \xi_0| - (x - x_0) \cdot (\xi - \xi_0)\Big) \ d\xi dx d\xi_{0} dx_0 \\ < \iiiint\limits_{{\rm I\!R}^{4n}}f(x+t\xi, \xi, t)f(x_{0}+t\xi_0, \xi_{0}, t)  \Big(|x - x_0|^2+|\xi - \xi_0|^2 \\ - (x-x_0).(\xi - \xi_0)\Big)d\xi dx d\xi_0 dx_0
\end{multline*}
We will use the structure of the equation in (\ref{odeform}) and integrate the interactions along the characteristics. The estimate above leads to the following inequality: 
\begin{multline*}
 \iiiint\limits_{{\rm I\!R}^{4n}}f(x, \xi, t)f(x_{0}, \xi_{0}, t) \Big(|x - x_0||\xi - \xi_0| - (x - x_0) \cdot (\xi - \xi_0)\Big) \ d\xi dx d\xi_{0} dx_0  \\ < \iiiint\limits_{{\rm I\!R}^{4n}}f(x, \xi, 0)f(x_{0}, \xi_{0}, 0) \Big(|x - x_0|^2+|\xi - \xi_0|^2 \\ - (x-x_0).(\xi - \xi_0)\Big)  d\xi dx d\xi_0 dx_0 + Z_1 + Z_2
 \end{multline*}
 In the expression above $Z_1$ and $Z_2$ are defined as: 
 \begin{multline*}
 Z_1 = 2\iiiint\limits_{{\rm I\!R}^{4n}}\int_{0}^{t} f(x, \xi, 0)I(x_{0}+s\xi_0, \xi_{0}, s) \Big(|x - x_0|^2+|\xi - \xi_0|^2 \\ - (x-x_0).(\xi - \xi_0)\Big) ds d\xi dx d\xi_0 dx_0
 \end{multline*}
 \begin{multline*}
 Z_2=   \iiiint\limits_{{\rm I\!R}^{4n}}\int_{0}^{t}\int_{0}^{t} I(x+s\xi, \xi, s)I(x_{0}+z\xi_0, \xi_{0}, z) \Big(|x - x_0|^2+|\xi - \xi_0|^2 \\- (x-x_0).(\xi - \xi_0)\Big)ds dz d\xi dx d\xi_0 dx_0
\end{multline*}
Using Definition \ref{mesoscopic}, we conclude that $Z_1$ and $Z_2$ are zero. Therefore we complete the proof of the first part of the theorem: 
\begin{multline*}
     \|f\|_{IG} < \iiiint\limits_{{\rm I\!R}^{4n}}f(x, \xi, 0)f(x_{0}, \xi_{0}, 0) \Big(|x - x_0|^2+|\xi - \xi_0|^2 \\- (x - x_0) \cdot (\xi - \xi_0)\Big) \ d\xi dx d\xi_{0} dx_0 
\end{multline*}
We just proved that $\|f\|_{IG}$ is bounded, therefore as a consequence of Lemma \ref{lemmaunc} and Definition \ref{ian}, we conclude that interaction uncertainty goes to infinity with time and complete the proof of the theorem: 
\begin{align*}
\lim_{t\rightarrow\infty}U_{I}(t) = \infty
\end{align*}
\end{proof}
\end{theorem}

\newpage
\section{Concentration of interactions} 
In this section we will use the bound we found for $\|f\|_{IG}$ from Section 3 along with the kinetic interaction Morawetz estimate (Theorem \ref{kime}) to demonstrate a result about the concentration of 2-particle interactions over time. \\

\indent Recall the notion of a blind cone defined in (\ref{bcone}). We will define a punctured blind cone by including the velocities that the magnitude of differences between them and some fixed $\xi_0$ are less than some positive constant: 
\begin{definition}\label{pbcone}
Let $C_{x_{0}, \xi_0}(x,c)$ be as defined in (\ref{bcone}). Define $K_{x_0, \xi_0}(x, c) \subset {\rm I\!R}^{n}$ as the {\it punctured blind cone} at point $x$ with respect to the observer $x_0$ that moves with the velocity $\xi_0$ for a fixed positive $v$:
\begin{align*}
    K_{x_{0},\xi_0}(x, c, v) = C_{x_{0}, \xi_0}(x,c) \cup B(\xi_0, v)    
\end{align*}
\end{definition}

\begin{definition}\label{gamma}
Let $\Gamma_{x_{0}, \xi_0}(c, v) \subset {\rm I\!R}^{4n}$ be the {\it collection of punctured blind cones} over the spatial variable  for some positive $c$ and $v$.
\begin{align*}
    \Gamma(c,v)=\{(x_0, \xi_0, x,\xi) \in {\rm I\!R}^n \times {\rm I\!R}^n \times {\rm I\!R}^n \times {\rm I\!R}^n | \ \xi \in K_{x_{0}, \xi_0}(x, c, v)\} 
\end{align*}
\end{definition}
\begin{remark}
An illustration of $\Gamma(c,v)$ for the specific case of an idle observer ($\xi_0=0$) at the origin ($x_0 = 0$) can be found in \cite{moini2021evolution}. \\
\end{remark}

\indent As a consequence of the conservation of mass, we have that the the total mass of 2-particle interactions remains constant overtime, therefore the following average is true: 
\begin{equation}\label{m^2}
\begin{aligned}
    \lim_{T\rightarrow\infty} \frac{1}{T}\int_{0}^{T}\iiiint\limits_{{\rm I\!R}^{4n}}f(x, \xi, t)f(x_{0}, \xi_{0}, t) d\xi dx d\xi_{0} d\xi_0 dt = {\bf M^2}
\end{aligned}
\end{equation}
In the theorem below, we will show that these interactions are actually being concentrated within the collection of punctured blind cones with arbitrary small $c$ and $v$. $\Gamma(c,v)\subset {\rm I\!R}^{4n}$ is much smaller than ${\rm I\!R}^{4n}$ and as the constants $c$ and $v$ go to zero, it converges to a measure zero subset. 

\begin{theorem}\label{conc} [{\it Concentration of 2-particle interactions}] Assume $f$ is a positive solution of the equation (1.1) and that the assumptions (\ref{mme}) and (\ref{localization}) are true. The 2-particle interactions, averaged over time, will almost always occur within the collection of punctured blind cones $\Gamma(c, v)$ for any positive $c$ and $v$:
\begin{align*}
\lim_{T\rightarrow\infty} \frac{1}{T}\int_{0}^{T}\iiiint\limits_{\Gamma(c,v)}f(x, \xi, t)f(x_{0}, \xi_{0}, t) d\xi dx d\xi_{0} d\xi_0 dt={\bf M^2}
    \end{align*}
\begin{proof}
Consider the following two subsets of ${\rm I\!R}^{4n}$: 
\begin{align*}
    N_R(c,v) =\{ (x_0, \xi_0, x, \xi)\in {\rm I\!R}^n \times {\rm I\!R}^n \times {\rm I\!R}^n \times ({\rm I\!R}^n - K_{x_0, \xi_0}(x,c,v)) \ |\ \  |x - x_0|>R\} \\ 
  M_R(c,v) =\{ (x_0, \xi_0, x, \xi)\in {\rm I\!R}^n \times {\rm I\!R}^n \times {\rm I\!R}^n \times ({\rm I\!R}^n - K_{x_0, \xi_0}(x,c,v)) \ |\ \  |x - x_0|\leq R\}
\end{align*}
It is possible to divide ${\rm I\!R}^{4n}$ into 3 mutuality exclusive sets $N_R, M_R$ and $\Gamma(c,v)$. Therefore the equation (\ref{m^2}) leads to: 
\begin{multline*}
      \lim_{T\rightarrow\infty} \frac{1}{T}\int_{0}^{T}\iiiint\limits_{{\rm I\!R}^{4n}}f(x, \xi, t)f(x_{0}, \xi_{0}, t) d\xi dx d\xi_{0} d\xi_0 dt \\=\lim_{T\rightarrow\infty} \frac{1}{T}\int_{0}^{T}\Big( \iiiint\limits_{\Gamma(c,v)}f(x, \xi, t)f(x_{0}, \xi_{0}, t) dx d\xi dx_{0} d\xi_{0} \\ + \iiiint\limits_{M_R(c,v)} f(x, \xi, t)f(x_{0}, \xi_{0}, t)dx d\xi dx_{0} d\xi_{0} \\+ \iiiint\limits_{N_R(c,v)}f(x, \xi, t)f(x_{0}, \xi_{0}, t)dx d\xi dx_{0} d\xi_{0} \Big) dt = {\bf M^2}
\end{multline*}
 The kinetic interaction Morawetz estimate (Theorem \ref{kime}) implies the estimate below for the interactions within $M_R(c,v)$:

\begin{align*}
    \int_{0}^{\infty}\iiiint\limits_{M_R(c,v)} f(x, \xi, t)f(x_{0}, \xi_{0}, t)dx d\xi dx_{0} d\xi_{0}dt < \infty
\end{align*}
Consequently we get: 
\begin{multline*}
          \lim_{T\rightarrow\infty} \frac{1}{T}\int_{0}^{T}\iiiint\limits_{{\rm I\!R}^{4n}}f(x, \xi, t)f(x_{0}, \xi_{0}, t) d\xi dx d\xi_{0} d\xi_0 dt \\=\lim_{T\rightarrow\infty} \frac{1}{T}\int_{0}^{T}\Big( \iiiint\limits_{\Gamma(c,v)}f(x, \xi, t)f(x_{0}, \xi_{0}, t) dx d\xi dx_{0} d\xi_{0} \\ + \iiiint\limits_{N_R(c,v)}f(x, \xi, t)f(x_{0}, \xi_{0}, t)dx d\xi dx_{0} d\xi_{0} \Big) dt = {\bf M^2}
\end{multline*}

\noindent It is possible to use the interaction relative angular norm (Definition \ref{ian} and Theorem \ref{ian}) to put an upper bound for the interactions within $N_R(c,v)$: 
\begin{multline*}
    (Rv - \cos(c))\iiiint\limits_{N_R(c,v)}f(x, \xi, t)f(x_{0}, \xi_{0}, t)dx d\xi dx_{0} d\xi_{0}  \\ < \iiiint\limits_{N_R(c,v)}f(x, \xi, t)f(x_{0}, \xi_{0}, t) \Big(|x - x_0||\xi - \xi_0| - (x - x_0) \cdot (\xi - \xi_0 \Big)dx d\xi dx_{0} d\xi_{0}
\end{multline*}
Therefore:
\begin{align*}
    \iiiint\limits_{N_R(c,v)}f(x, \xi, t)f(x_{0}, \xi_{0}, t)dx d\xi dx_{0} d\xi_{0} < \frac{\|f\|_{IG}}{Rv - \cos(c)}
\end{align*}
The inequality above and the previous computations imply: 
\begin{align*}
{\bf M^2} - \frac{\|f\|_{IG}}{Rv(1-\cos(c))}  \leq \liminf_{T\rightarrow\infty} \frac{1}{T}\int_{0}^{T}(\iiiint\limits_{\Gamma(c,v)} f(x,\xi, t)\ dxd\xi)dt \leq {\bf M}^2
\end{align*}
Since the former inequality  is valid for any arbitrary $R>0$, we complete the proof:
\begin{align*}
\lim_{T\rightarrow\infty} \frac{1}{T}\int_{0}^{T}\iiiint\limits_{\Gamma(c,v)}f(x, \xi, t)f(x_{0}, \xi_{0}, t) d\xi dx d\xi_{0} d\xi_0 dt={\bf M^2}
\end{align*}
\end{proof}
\end{theorem}

\printbibliography
\end{document}